\documentclass[12pt]{article}
\usepackage{graphicx}

\usepackage{epsfig}
\usepackage{amsmath}
\usepackage{amssymb}
\usepackage{amsthm}
\usepackage{latexsym}
\usepackage{cite}
\usepackage{amsbsy,amsfonts,euscript,eepic,rotating}

\newcommand{\comment}[1]{}
\newtheorem{theorem}{Theorem}
\newtheorem{lemma}{Lemma}[section]

\newtheorem{definition}{Definition}[section]
\begin{document}

\title{\LARGE
{\bf $\text{SLE}_6$ and $\text{CLE}_6$ from Critical Percolation}
\thanks{Research supported in part by the NSF under grant PHY99-07949
(preprint no.~NSF-KITP-06-76).}}

\author{
{\bf Federico Camia}
\thanks{Research supported in part by a Veni grant of the
NWO (Dutch
Organization for Scientific Research).}\,
\thanks{E-mail: fede@few.vu.nl}\\
{\small \sl Department of Mathematics, Vrije Universiteit Amsterdam}\\
\and
{\bf Charles M.~Newman}
\thanks{Research supported in part by the NSF under grants
DMS-01-04278 and DMS-06-06696.}\,
\thanks{E-mail: newman@courant.nyu.edu}\\
{\small \sl Courant Inst.~of Mathematical Sciences,
New York University}
}

\date{}

\maketitle

\begin{abstract}
We review some of the recent progress on the scaling
limit of two-dimensional critical percolation; in particular,
the convergence of the exploration path to chordal SLE$_6$
and the ``full" scaling limit of cluster interface loops.
The results given here on the full scaling limit and its
conformal invariance extend those presented previously.
For site percolation on the triangular lattice, the results
are fully rigorous.
We explain some of the main ideas, skipping most technical details.
\end{abstract}



\section{Introduction} \label{intro}
In the theory of critical phenomena it is usually assumed
that a physical system near a continuous phase transition
is characterized by a single length scale (the ``correlation
length'') in terms of which all other lengths should be measured.
When combined with the experimental observation that the
correlation length diverges at the phase transition, this
simple but strong assumption, known as the scaling hypothesis,
leads to the belief that at criticality the system has
no characteristic length, and is therefore invariant under
scale transformations.
This suggests that all thermodynamic functions at criticality
are homogeneous functions, and predicts the appearance of
power laws.

It also implies that if one rescales appropriately a critical lattice
model, shrinking the lattice
spacing to zero, it should be possible to obtain
a continuum model, known as the ``scaling limit."
The scaling limit is not restricted to a lattice and may possess more
symmetries than the original model.
Indeed, the scaling limits of many critical lattice models are believed
to be conformally invariant and to correspond to Conformal Field Theories
(CFTs).
But until recently, such a correspondence was at most heuristic, and was
assumed as a starting point by physicists working in CFT.
The methods of CFT themselves proved hard to put into a rigorous mathematical
formulation.

The introduction by Oded Schramm~\cite{schramm} of the Stochastic/Schramm
Loewner Evolution (SLE) has provided a new powerful and
mathematically rigorous tool to study scaling limits of critical lattice models.
Thanks to this, in recent years tremendous progress has been made in understanding
the conformally invariant nature of the scaling limits of several such models.

While CFT focuses on correlation functions of local ``operators" (e.g., spin
variables in the Ising model), SLE describes the behavior of macroscopic random
curves present in these models, such as percolation cluster boundaries.
In the scaling limit, the distribution of such random curves can be uniquely
identified thanks to their conformal invariance and a certain ``Markovian" property.
There is a one-parameter family of SLEs, indexed by a positive real number
$\kappa$, and they appear to be essentially the only possible candidates for
the scaling limits of interfaces of two-dimensional critical systems that are
believed to be conformally invariant.

The main power of SLE stems from the fact that it allows to compute different
quantities; for example, percolation crossing probabilities and various
percolation critical exponents.
Therefore, relating the scaling limit of a critical lattice model to SLE allows
for a rigorous determination of some aspects of the large scale behavior of the
lattice model.

In the context of the Ising, Potts and $O(n)$ models, an SLE curve is believed
to describe the scaling limit of a single interface, which can be obtained by
imposing special boundary conditions.
A single SLE curve is therefore not in itself sufficient to immediately describe
the scaling limit of the unconstrained model without boundary conditions in the
whole plane (or in domains with boundary conditions that do not determine a single
interface), and contains only limited information concerning the connectivity
properties of the model.

A more complete description can be obtained in terms of loops, corresponding
to the scaling limits of cluster boundaries.
Such loops should also be random and have a conformally invariant distribution.
This approach led Wendelin Werner~\cite{werner5,werner6}
(see also~\cite{werner3}) to the definition
of Conformal Loop Ensembles (CLEs), which are, roughly speaking, random
collections of fractal loops with a certain
``conformal restriction/renewal property."

For percolation, a complete proof of the connection with SLE, first conjectured
by Schramm in~\cite{schramm}, has recently been given in~\cite{cn2}.
The proof relies heavily on the ground breaking result of Stas
Smirnov~\cite{smirnov,smirnov-long} about the existence and conformal invariance
of the scaling limit of crossing probabilities (see~\cite{cardy1}).
The last section of this paper explains the main ideas of that proof,
highlighting the role of conformal invariance, but without dwelling on
the heavy technical details.

As for the Ising, Potts and $O(n)$ models, the scaling limit of percolation in
the whole plane should be described by a measure on loops, where the loops are
closely related to SLE curves.
Such a description in the case of percolation was presented in~\cite{cn},
where the authors of the present paper constructed a probability measure
on collections of fractal conformally invariant loops in the plane (closely
related to a CLE), arguing that it corresponds to the ``full" scaling limit
of critical two-dimensional percolation.
A proof of that statement was subsequently provided in~\cite{cn1}.

Here, we will briefly explain how to go from a single SLE curve to the full
scaling limit, again skipping the technical details,
for the case of a Jordan domain with ``monochromatic" boundary conditions
(see Theorem~\ref{thm-full-scaling-limit}).
This extends the results presented in~\cite{cn1}, where the scaling limit
was first taken in the unit disc and then an infinite volume limit was
taken in order to obtain the full scaling limit in the whole plane.
Moving from the unit disc (or any convex domain) to a general Jordan
domain introduces extra complications that are dealt with using a new
argument, developed in~\cite{cn2}, that exploits the continuity of Cardy's
formula~\cite{cardy1} with respect to changes in the shape of the domain
(see the discussion in Sec.~\ref{sec-convergence-full}).
Taking scaling limits in general Jordan domains is a necessary step
in order to consider ``conformal restriction/renewal properties" as
in Theorem~\ref{thm-cle} below.

Using the full scaling limit, one can attempt to understand the geometry of
the ``near-critical" scaling limit, where the percolation density tends to
the critical one in an appropriate way as the lattice spacing tends to zero.
A heuristic analysis~\cite{cfn1,cfn2} based on a natural ansatz leads to a
one-parameter family of loop models (i.e., probability measures on
random collections of loops), with the critical full scaling limit
corresponding to a particular choice of the parameter.
Except for the latter case, these measures are not scale invariant, but are
mapped into one another by scale transformations.
This framework can be used to define a renormalization group flow (under the
action of dilations), and to describe the scaling limit of related models,
such as invasion and dynamical percolation and the minimal spanning tree.
In particular, 
this analysis helps explain why the scaling limit of the minimal spanning
tree may be scale invariant but \emph{not} conformally invariant, as first
observed numerically by Wilson~\cite{wi}.

\section{SLE and CLE} \label{sec-sle-cle}


The Stochastic/Schramm Loewner Evolution with parameter $\kappa > 0$
($\text{SLE}_{\kappa}$) was introduced by Schramm~\cite{schramm} as
a tool for studying the scaling limit of two-dimensional discrete
(defined on a lattice) probabilistic models whose scaling limits are
expected to be conformally invariant.
In this section we define the chordal version of $\text{SLE}_{\kappa}$;
for more on the subject, the interested reader can consult
the original paper~\cite{schramm} as well as the fine
reviews by Lawler~\cite{lawler1}, Kager and Nienhuis~\cite{kn},
and Werner~\cite{werner4}, and Lawler's book~\cite{lawler2}.

Let $\mathbb H$ denote the upper half-plane.
For a given continuous real function $U_t$ with $U_0 = 0$,
define, for each $z \in \overline{\mathbb H}$, the function
$g_t(z)$ as the solution to the ODE
\begin{equation}
\partial_t g_t(z) = \frac{2}{g_t(z) - U_t},
\end{equation}
with $g_0(z) = z$.
This is well defined as long as $g_t(z) - U_t \neq 0$,
i.e., for all $t < T(z)$, where
\begin{equation}
T(z) \equiv \sup \{ t \geq 0 : \min_{s \in [0,t]} | g_s(z) - U_s| > 0 \}.
\end{equation}
Let $K_t \equiv \{ z \in \overline{\mathbb H} : T(z) \leq t \}$
and let ${\mathbb H}_t$ be the unbounded component of
${\mathbb H} \setminus K_t$; it can be shown that $K_t$ is bounded
and that $g_t$ is a conformal map from ${\mathbb H}_t$ onto $\mathbb H$.
For each $t$, it is possible to write $g_t(z)$ as
\begin{equation}
g_t(z) = z + \frac{2t}{z} + O\left(\frac{1}{z^2}\right),
\end{equation}
when $z \to \infty$.
The family $(K_t, t \geq 0)$ is called the {\bf Loewner chain}
associated to the driving function $(U_t, t \geq 0)$.

\begin{definition} \label{def-sle}
{\bf Chordal $\text{SLE}_{\kappa}$} is the Loewner chain $(K_t, t \geq 0)$
that is obtained when the driving function
$U_t = \sqrt{\kappa} B_t$ is $\sqrt{\kappa}$ times a standard
real-valued Brownian motion $(B_t, t \geq 0)$ with $B_0 = 0$.
\end{definition}

For all $\kappa \geq 0$, chordal $\text{SLE}_{\kappa}$ is almost surely generated
by a continuous random curve $\gamma$ in the sense that, for all $t \geq 0$,
${\mathbb H}_t \equiv {\mathbb H} \setminus K_t$ is the unbounded connected
component of ${\mathbb H} \setminus \gamma[0,t]$; $\gamma$ is called the
{\bf trace} of chordal $\text{SLE}_{\kappa}$.

It is not hard to see, as argued by Schramm, that any continuous random curve
$\gamma$ in the upper half-plane starting at the origin and going to infinity
must be an SLE curve if it possesses the following ``{\bf conformal Markov property}."
For any fixed $T \in {\mathbb R}$, conditioning on $\gamma[0,T]$, the image
under $g_T$ of $\gamma[T,\infty)$ is distributed like an independent copy of
$\gamma$, up to a time reparametrization.
This implies that the driving function $U_t$ in the Loewner chain associated to
the curve $\gamma$ is continuous and has stationary and independent increments.
If the time parametrization implicit in Definition~\ref{def-sle} and the discussion
preceding it is chosen for $\gamma$, then scale invariance also implies that the
law of $U_t$ is the same as the law of $\lambda^{-1/2} U_{\lambda t}$ when $\lambda>0$.
These properties together imply that $U_t$ must be a constant multiple of standard
Brownian motion.


Let now $D \subset {\mathbb C}$ ($D \neq {\mathbb C}$) be a simply
connected domain whose boundary is a continuous curve.
By Riemann's mapping theorem, there are (many) conformal maps
from the upper half-plane $\mathbb H$ onto $D$.
In particular, given two distinct points $a,b \in \partial D$
(or more accurately, two distinct prime ends), there exists a
conformal map $f$ from $\mathbb H$ onto $D$ such that $f(0)=a$
and $f(\infty) \equiv \lim_{|z| \to \infty} f(z) = b$.
In fact, the choice of the points $a$ and $b$ on the boundary
of $D$ only characterizes $f(\cdot)$ up to a scaling factor $\lambda>0$,
since $f(\lambda \cdot)$ would also do.

Suppose that $(K_t, t \geq 0)$ is a chordal $\text{SLE}_{\kappa}$ in
$\mathbb H$ as defined above; we define chordal $\text{SLE}_{\kappa}$
$(\tilde K_t, t \geq 0)$ in $D$ from $a$ to $b$ as the
image of the Loewner chain $(K_t, t \geq 0)$ under $f$.
It is possible to show, using scaling properties of
$\text{SLE}_{\kappa}$, that the law of $(\tilde K_t, t \geq 0)$
is unchanged, up to a linear time-change, if we replace
$f(\cdot)$ by $f(\lambda \cdot)$.
This makes it natural to consider $(\tilde K_t, t \geq 0)$ as
a process from $a$ to $b$ in $D$, ignoring the role of $f$.
The trace of chordal SLE in $D$ from $a$
to $b$ will be denoted by $\gamma_{D,a,b}$.


We now move from the conformally invariant random curves of SLE to
collections of conformally invariant random loops and introduce the concept
of Conformal Loop Ensemble (CLE --- see~\cite{werner3,werner5,werner6,sheffield}).
The key feature of a CLE is a sort of ``conformal restriction/renewal property."
Roughly speaking, a CLE in $D$ is a random collection ${\cal L}_D$ of loops
such that if all the loops intersecting a (closed) subset $D'$ of $D$ or of
its boundary are removed, the loops in any one of the various remaining
(disjoint) subdomains of $D$ form a random collection of loops distributed
as an independent copy of ${\cal L}_D$ conformally mapped to that subdomain
(see Theorem~\ref{thm-cle}).
We will not attempt to be more precise here since somewhat different definitions
(although, in the end, substantially equivalent) have appeared in the literature,
but the meaning of the conformal restriction/renewal property should be clear
from Theorem~\ref{thm-cle}.

For formal definitions and more discussion on the properties of a CLE,
the reader is referred to the original literature on the subject~\cite{werner5,werner6,sheffield},
where it is shown that there is a one-parameter family $\text{CLE}_{\kappa}$
of conformal loop ensembles with the above conformal restriction/renewal
property and that for $\kappa \in (8/3,8]$, the $\text{CLE}_{\kappa}$ loops
locally ``look like" $\text{SLE}_{\kappa}$ curves.

There are numerous lattice models that can be described in terms of random
curves and whose scaling limits are assumed (and in a few cases proved) to
be conformally invariant.
These include the Loop Erased Random Walk, the Self-Avoiding Walk and the
Harmonic Explorer, all of which can be defined as polygonal paths along
the edges of a lattice.
The Ising, Potts and percolation models instead are naturally defined in
terms of clusters, and the interfaces between different clusters form random
loops.
In the $O(n)$ model, configurations of loops along the edges of the hexagonal
lattice are weighted according to the total number and length of the loops.
All of these models are supposed to have scaling limits described by
$\text{SLE}_{\kappa}$ or $\text{CLE}_{\kappa}$ for some value of $\kappa$
between $2$ and $8$.
For more information on these lattice models and their scaling limits, the
interested reader can consult~\cite{kn,cardy2,cardy3,werner5,sheffield}.

In the
rest of the paper we will restrict attention to percolation, where the
connection with $\text{SLE}_6$ and $\text{CLE}_6$ has been made rigorous
\cite{smirnov,smirnov-long,cn1,cn2}.

\section{Conformal Invariance of Critical Percolation} \label{sec-conformal}


In this section we will consider critical site percolation on the
triangular lattice, for which conformal invariance in the scaling limit was
rigorously proved~\cite{smirnov,smirnov-long}.
A precise formulation of conformal invariance, attributed to Michael Aizenman,
is that the probability that a percolation cluster crosses between two disjoint
segments of the boundary of some simply connected domain should converge to a
conformally invariant function of the domain and the two segments of the boundary.
This conjecture is connected with the extensive numerical investigations
reported in~\cite{lps}.
A formula for the purposed limit was then derived~\cite{cardy1} by John Cardy
using
(non-rigorous) field theoretical methods.
The interest of mathematicians was already evident in~\cite{lps}, but a proof
of the conjecture~\cite{smirnov,smirnov-long} (and of Cardy's formula) did not
come until 2001.

\begin{figure}[!ht]
\begin{center}
\includegraphics[width=6cm]{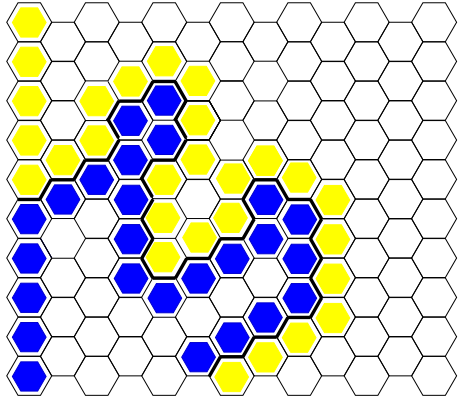}
\caption{Percolation exploration path in a portion of the
hexagonal lattice with blue/yellow boundary conditions on the
first column, corresponding to the boundary of the region where
the exploration is carried out.
The colored hexagons that do not belong to the first column
have been ``explored'' during the exploration process.
The heavy line between yellow (light) and blue (dark) hexagons
is the exploration path produced by the exploration process.}
\label{fig-exp-path}
\end{center}
\end{figure}

We will denote by $\cal T$ the two-dimensional triangular lattice,
whose sites are identified with the elementary cells of a regular hexagonal
lattice $\cal H$ embedded in the plane as in Fig.~\ref{fig-exp-path}.
We say that two hexagons are neighbors (or that they are adjacent) if
they have a common edge.
A sequence $(\xi_0, \ldots, \xi_n)$ of 
hexagons of $\cal H$ such that $\xi_{i-1}$ and $\xi_i$
are neighbors for all $i= 1, \ldots, n$ and $\xi_i \neq \xi_j$
whenever $i \neq j$ will be called a {\bf $\cal T$-path}.
If the first and last hexagons of the path are neighbors,
the path will be called a {\bf $\cal T$-loop}.

Let $D$ be a bounded simply connected domain containing the origin
whose boundary $\partial D$ is a continuous curve.
Let $\phi:\overline{\mathbb D} \to D$ be the (unique) continuous function
that maps $\mathbb D$ onto $D$ conformally and such that $\phi(0)=0$ and
$\phi'(0)>0$.
Let $z_1,z_2,z_3,z_4$ be four points of $\partial D$ in counterclockwise
order --- i.e., such that $z_j=\phi(w_j), \,\,\, j=1,2,3,4$, with
$w_1,\ldots,w_4$ in counterclockwise order.
Also, let $\eta = \frac{(w_1-w_2)(w_3-w_4)}{(w_1-w_3)(w_2-w_4)}$.
Cardy's formula~\cite{cardy1} for the probability $\Phi_{D}(z_1,z_2;z_3,z_4)$
of a ``crossing" inside $D$ from the counterclockwise arc $\overline{z_1 z_2}$
to the counterclockwise arc $\overline{z_3 z_4}$ is
\begin{equation} \label{cardy-formula}
\Phi_{D}(z_1,z_2;z_3,z_4) =
\frac{\Gamma(2/3)}{\Gamma(4/3) \Gamma(1/3)} \eta^{1/3} {}_2F_1(1/3,2/3;4/3;\eta),
\end{equation}
where ${}_2F_1$ is a hypergeometric function.

For a given mesh $\delta>0$, the probability of a blue crossing inside $D$
from the counterclockwise arc $\overline{z_1 z_2}$ to the counterclockwise
arc $\overline{z_3 z_4}$ is the probability of the existence of a blue
$\cal T$-path $(\xi_0,\ldots,\xi_n)$ such that $\xi_0$ intersects the
counterclockwise arc $\overline{z_1 z_2}$, $\xi_n$ intersects the
counterclockwise arc $\overline{z_3 z_4}$, and $\xi_1,\ldots,\xi_{n-1}$
are all contained in $D$.
Smirnov~\cite{smirnov,smirnov-long} proved that crossing probabilities
converge in the scaling limit to conformally invariant functions of the
domain and the four points on its boundary, and identified the limit
with Cardy's formula~(\ref{cardy-formula}).



The proof of Smirnov's theorem is based on the identification of certain
generalized crossing probabilities that are almost discrete harmonic functions
and whose scaling limits converge to harmonic functions.
The behavior on the boundary of such functions is easy to determine and is
sufficient to specify them uniquely.
The relevant crossing probabilities can be expressed in terms of the boundary
values of such harmonic functions, and as a consequence are invariant under
conformal transformations of the domain and the two segments of its boundary.

The presence of a blue crossing in $D$ from the counterclockwise boundary arc
$\overline{z_1 z_2}$ to the counterclockwise boundary arc $\overline{z_3 z_4}$
can be determined using a clever algorithm that explores the percolation
configuration inside $D$ starting at, say, $z_1$ and assumes that the hexagons
just outside $\overline{z_1 z_2}$ are all blue and those just outside
$\overline{z_4 z_1}$ are all yellow.
The exploration proceeds following the interface between the blue cluster
adjacent to $\overline{z_1 z_2}$ and the yellow cluster adjacent to
$\overline{z_4 z_1}$.
A blue crossing is present if the exploration process reaches $\overline{z_3 z_4}$
before $\overline{z_2 z_3}$.
This {\bf exploration process} and the {\bf exploration path}
(see Fig.~\ref{fig-exp-path}) associated to it were introduced by
Schramm in~\cite{schramm}.

The exploration process can be carried out in ${\mathbb H} \cap {\cal H}$,
where the hexagons in the lowest row and to the left of a chosen hexagon
have been colored yellow and the remaining hexagons in the lowest row
have been colored blue.
This produces an infinite exploration path, whose scaling limit was
conjectured~\cite{schramm} by Schramm to converge to $\text{SLE}_6$.

It is easy to see that the exploration process is Markovian in the sense
that, conditioned on the exploration up to a certain (stopping) time,
the future of the exploration evolves in the same way as the past except
that it is now performed in a different domain, where some of the the
explored hexagons have become part of the boundary (see, e.g., Fig.~\ref{fig-exp-path}).

This observation, together with the connection between the exploration
process and crossing probabilities, Smirnov's theorem about the conformal
invariance of crossing probabilities in the scaling limit, and Schramm's
characterization of SLE via the conformal Markov property discussed in
Sec.~\ref{sec-sle-cle}, strongly support the above conjecture.

As we now explain, the natural setting to define the exploration process
is that of {\bf lattice domains}, i.e., sets $D^{\delta}$ of hexagons of
$\delta\cal H$ that are {\bf connected} in the sense that any two hexagons
in $D^{\delta}$ can be joined by a
$(\delta\cal T)$-path contained in $D^{\delta}$.
We say that a bounded lattice domain $D^{\delta}$ is {\bf simply connected}
if both $D^{\delta}$ and $\delta{\cal T} \setminus D^{\delta}$ are connected.
A {\bf lattice-Jordan} domain $D^{\delta}$ is a bounded simply connected
lattice domain such that the set of hexagons adjacent to $D^{\delta}$
is a $(\delta\cal T)$-loop.

Given a lattice-Jordan domain $D^{\delta}$, the set of hexagons
adjacent to $D^{\delta}$ can be partitioned into two (lattice-)connected sets.
If those two sets of hexagons are assigned different colors,
for any coloring of the hexagons inside $D^{\delta}$, there is an
interface between two clusters of different colors starting and
ending at two boundary points, $a^{\delta}$ and $b^{\delta}$,
corresponding to the locations on the boundary of $D^{\delta}$
where the color changes.
If one performs an exploration process in $D^{\delta}$ starting
at $a^{\delta}$, one ends at $b^{\delta}$, producing an exploration
path $\gamma^{\delta}$ that traces the entire interface from
$a^{\delta}$ to $b^{\delta}$.

Given a planar domain $D$, we denote by $\partial D$ its topological boundary.
Let $\partial D$ be locally connected (i.e., a continuous curve),
and assume that $D$ contains the origin.
Then one can parametrize $\partial D$ by $\varphi:S^1 \to \partial D$,
where $\varphi$ is the restriction to the unit circle $S^1$ of the
continuous map $\phi:\overline{\mathbb D} \to \overline D$ that
is conformal in $\mathbb D$ and satisfies $\phi(0)=0$, $\phi'(0)>0$.
With this notation, we say that $D^{\delta}$ converges to $D$ as
$\delta \to 0$ if
\begin{equation} \label{domain-convergence}
\lim_{\delta \to 0}\inf_h \sup_{z \in S^1} |\varphi(z)-\varphi^{\delta}(h(z))| =0,
\end{equation}
where the infimum is over monotonic functions $h:S^1 \to S^1$ (and
the objects with the superscript $\delta$ refer to $D^{\delta}$ ---
for simplicity we are assuming that all domains contain the origin).
If moreover two points, $a^{\delta},b^{\delta} \in \partial D^{\delta}$,
converge respectively to $a,b \in \partial D$ as $\delta \to 0$,
we write $(D^{\delta},a^{\delta},b^{\delta}) \to (D,a,b)$.
In the following theorem the topology on curves is that induced by
the supremum norm, but with monotonic reparametrizations of the curves
allowed (see~\cite{ab,cn1,cn2}), i.e., the distance between curves is
\begin{equation} \label{distance}
\text{d}(\gamma,\gamma^{\delta})=\inf_h \sup_{t \in [0,\infty)} |\gamma(t)-\gamma^{\delta}(h(t))|,
\end{equation}
where $\gamma(t),\gamma^{\delta}(t), t \in [0,\infty)$, are parametrizations
of $\gamma_{D,a,b}$ and $\gamma^{\delta}_{D,a,b}$ respectively,
and the infimum is over monotonic functions $h:[0,\infty) \to [0,\infty)$.
A proof of the theorem can be found in~\cite{cn2} and a detailed
sketch is presented in Sec.~\ref{sec-convergence-sle} below.

\begin{theorem} \label{thm-conv-to-sle}
Let $(D,a,b)$ be a Jordan domain with two distinct selected points
on its boundary $\partial D$.
Then, for lattice-Jordan domains $D^{\delta}$ from $\delta{\cal H}$ with
$a^{\delta},b^{\delta} \in \partial D^{\delta}$ such that
$(D^{\delta},a^{\delta},b^{\delta}) \to (D,a,b)$ as $\delta \to 0$,
the percolation exploration path $\gamma^{\delta}_{D,a,b}$ in $D^{\delta}$
from $a^{\delta}$ to $b^{\delta}$ converges in distribution to the
trace $\gamma_{D,a,b}$ of chordal $\text{\emph{SLE}}_6$ in $D$
from $a$ to $b$, as $\delta \to 0$.
\end{theorem}

\section{The Full Scaling Limit in a Jordan Domain} \label{sec-full}


In this section we define the {\bf Continuum Nonsimple Loop (CNL) process}
in a Jordan domain $D$, a random collection of countably many nonsimple
fractal loops in D which corresponds to the ``full" scaling limit of
percolation in $D$ with monochromatic boundary conditions.
This refers to the collection of all cluster boundaries of percolation
configurations in $D$ with the hexagons at the boundary of $D$ all
blue (obviously, one could as well choose yellow boundary conditions).
The algorithmic construction that we present below is analogous to
that of~\cite{cn,cn1} for the unit disc $\mathbb D$, but here we
perform it in a general Jordan domain.

The CNL process on the full plane can be obtained by taking a sequence
of domains $D$ tending to ${\mathbb C}$.
This was done in~\cite{cn, cn1} and for that purpose, discs of radius
$R$ with $R \to \infty$ suffice.
This full plane CNL process is the scaling limit of the collection of
all cluster boundaries in the full lattice (without boundary conditions).
In order to consider conformal restriction/renewal properties (as we do
in Theorem~\ref{thm-cle} below), one needs to consider the CNL process
in fairly general bounded domains $D$.
There are extra complications in taking the scaling limit when $D$ is
non-convex, as discussed in Sec.~\ref{sec-convergence-full}.

The basic ingredient in our algorithmic construction consists of a
chordal $\text{SLE}_6$ path between two points on the boundary of
a Jordan domain.
As we will explain soon, sometimes the two boundary points are
``naturally'' determined as a product of the construction itself,
and sometimes they are given as an input to the construction.
In the second case, there are various procedures which would yield
the ``correct'' distribution for the resulting
CNL process; one possibility is as follows.
Given a domain $D$, choose $a$ and $b$ so that, of all points in
$\partial D$, they have maximal $x$-distance or maximal $y$-distance,
whichever is greater.
It is important to stress that in the end, the CNL process will turn
out to be independent of the actual choice of boundary points, as is
evident in Theorem~\ref{thm-full-scaling-limit}.
(One caveat is that one should avoid ``malicious" choices of the
boundary points for which the entire original domain would not be
explored asymptotically.)

The first step of our construction is a chordal $\text{SLE}_6$,
$\gamma\equiv\gamma_{D,a,b}$, between two boundary points
$a,b \in \partial D$ chosen according to the above rule.
The set $D \setminus \gamma_{D,a,b}[0,\infty)$ is a countable union of its
connected components, which are open and simply connected.
If $z$ is a deterministic point in $D$, then with probability one, $z$ is
not touched by $\gamma$~\cite{rs} and so belongs to a unique one of these,
that we denote $D_{a,b}(z)$.
There are four kinds of components which may be usefully thought of in terms
of how a point $z$ in the interior of the component was first ``trapped'' at
some time $t_1$ by $\gamma[0,t_1]$ perhaps together with either the counterclockwise
arc $\partial_{a,b} D$ of $\partial D$ between $a$ and $b$ or the counterclockwise
arc $\partial_{b,a} D$ of $\partial D$ between $b$ and $a$: (1) those components
whose boundary contains a segment of $\partial_{b,a} D$ between two successive
visits at $\gamma(t_0)$ and $\gamma(t_1)$ to $\partial_{b,a} D$ (where here and
below $t_0<t_1$), (2) the analogous components with $\partial_{b,a} D$ replaced
by the other part of the boundary $\partial_{a,b} D$, (3) those components formed
when $\gamma(t_0) = \gamma(t_1)$ with $\gamma$ winding about $z$ in a counterclockwise direction between $t_0$ and $t_1$, and finally (4) the analogous clockwise components.

To conclude the first step, we consider all domains of type (1),
corresponding to excursions of the $\text{SLE}_6$ path from the portion
$\partial_{b,a} D$ of $\partial D$.
For each such domain $D'$, the points $a'$ and $b'$ on its boundary
are chosen to be respectively those points where the excursion ends
and where it begins, that is, for $D_{a,b}(z)$ we set $a'=\gamma((t_1(z))$
and $b'=\gamma(t_0(z))$.
We then ``run" chordal $\text{SLE}_6$ from $a'$ to $b'$.
The loop obtained by pasting together the excursion from $b'$ to $a'$
followed by the new $\text{SLE}_6$ path from $a'$ to $b'$ is one of our continuum
loops.
At the end of the first step, then, the procedure has generated countably
many loops that touch
$\partial_{b,a} D$; each of these loops touches
$\partial_{b,a} D$ but may or may not touch $\partial_{a,b} D$.

The last part of the first step also produces new domains, corresponding
to the connected components of $D' \setminus \gamma_{D',a',b'}[0,\infty)$
for all domains $D'$ of type (1).
Each one of these components, together with all the domains of type (2),
(3) and (4) previously generated, is to be used in the next step of the
construction, playing the role of the original domain $D$.
For each one of these domains, we choose the ``new $a$'' and ``new $b$''
on the boundary as explained before, and then continue with the construction.
Note that the ``new $a$'' and ``new $b$'' are chosen according to the
rule explained at the beginning of this section also for domains of
type (2), even though they are generated by excursions like the domains
of type (1).

This iterative procedure produces at each step a countable set of loops.
The limiting object, corresponding to the collection of all such loops,
is our basic process.
(Technically speaking, we should include also trivial loops fixed at
each $z \in D$ so that the collection of loops is closed in an appropriate
sense~\cite{ab}.)

As explained, the construction is carried out iteratively and can be
performed simultaneously on all the domains that are generated at each step.
We wish to emphasize, though, that the obvious monotonicity of the
procedure, where at each step new paths are added independently in different
domains, and new domains are formed from the existing ones, implies that
any other choice of the order in which the domains are used would give the
same result (i.e., produce the same limiting distribution), provided that
every domain that is formed during the construction is eventually used.

The main interest of the loop process defined above is in the following
theorem, where the topology on collections of loops is that of
Aizenman-Burchard~\cite{ab} (see also~\cite{cn1}).

\begin{theorem} \label{thm-full-scaling-limit}
In the scaling limit, $\delta \to 0$, the collection of all cluster boundaries
of critical site percolation on the triangular lattice in a Jordan domain
$D$ with monochromatic boundary conditions converges in distribution to the
Continuum Nonsimple Loop process in $D$.
\end{theorem}

A key property of the CNL process is conformal invariance.
\begin{theorem} \label{thm-conformal-invariance}
Let $D,D'$ be two Jordan domains and $f:\overline D \to {\overline D}'$
a continuous function that maps $D$ conformally onto $D'$.
Then the CNL process in $D'$ is distributed like the image under $f$
of the CNL process in $D$.
\end{theorem}

Moreover, as shown in the next theorem, the outermost loops of the CNL process
in a Jordan domain satisfy a conformal restriction/renewal property, as in the
definitions of the Conformal Loop Ensembles of Werner~\cite{werner5} and
Sheffield~\cite{sheffield}.
\begin{theorem} \label{thm-cle}
Let $D$ be a Jordan domain and ${\cal L}_D$ be the collection of
CN loops
in $\overline D$ that are not surrounded by any other loop.
Consider an arc $\Gamma$ of $\partial D$ and let ${\cal L}_{D,\Gamma}$ be
the set of loops of ${\cal L}_D$ that touch $\Gamma$.
Then, conditioned on ${\cal L}_{D,\Gamma}$, for any connected component
$D'$ of $D \setminus \overline{\cup \{ L:L \in {\cal L}_{D,\Gamma} \}}$,
the loops in $\overline{D'}$ form a random collection of loops distributed
as an independent copy of ${\cal L}_D$ conformally mapped to $D'$.
\end{theorem}

Yet another form of conformal invariance is illustrated by showing how to
obtain a (conformally invariant) $\text{SLE}_6$ curve from the CNL process.
Given a Jordan domain $D$ and two points $a,b \in \partial D$, let $\Gamma=\overline{ba}$
be the counterclockwise closed arc $\overline{ba}$ of $\partial D$.
Define ${\cal L}_D$ and ${\cal L}_{D,\Gamma}$ as in Theorem~\ref{thm-cle}.
For each $L \in {\cal L}_{D,\Gamma}$, going from $a$ to $b$ clockwise, there are
a first and a last point, $x$ and $y$ respectively, where $L$ intersects $\Gamma$.
We call the counterclockwise arc of $L$ between $x$ and $y$, $\overline{xy} (L)$,
a (counterclockwise) {\bf excursion} from $\overline{ba}$.
We call such an $\overline{xy} (L)$ a {\bf maximal} excursion if there is
no other excursion
from $\overline{ba}$ in (the closure of) the domain created by
$\overline{xy} (L)$ and the counterclockwise arc $\overline{yx}$ of $\partial D$.
The random curve obtained by ``pasting together" (in the order in which they are
encountered going from $a$ to $b$ clockwise) all such maximal excursions from $\overline{ba}$
is distributed like a chordal $\text{SLE}_6$ in $D$ from $a$ to $b$.

The procedure described above obviously requires some care, since there are countably
many such excursions and there is no such thing as the ``first" excursion encountered
from $a$ or the ``next" excursion.
What this means is that in order to properly define the curve, one needs to use a
limiting procedure.
Since it is quite obvious how to do it but rather tedious to explain, we leave the
details to the interested reader (see~\cite{cn1}).

\section{Convergence and Conformal Invariance of the Full Scaling Limit} \label{sec-convergence-full}


\noindent {\bf Sketch of the Proof of Theorem~\ref{thm-full-scaling-limit}.}
The first step in the proof of Theorem~\ref{thm-full-scaling-limit} is to note
that it follows directly from the work of Aizenman-Burchard~\cite{ab} that the
family of distributions of the collections of cluster boundaries in $D$ with
monochromatic boundary conditions is tight, as $\delta \to 0$, in the sense
of the induced Hausdorff metric on closed sets of curves based on the
metric~(\ref{distance}) for single curves (see~\cite{ab} and~\cite{cn1}), and
so there is convergence along subsequences $\delta_k \to 0$.
What needs to be proved is that the limiting distribution is that of the CNL
process, independently of the subsequence $\delta_k$.

The key to the proof is an algorithmic construction on the lattice which
parallels the continuum construction of Sec.~\ref{sec-full} used to define
the CNL process in $D$.
The construction takes place in a lattice-domain
$D_k \equiv D^{\delta_k}$ that
converges to $D$ in the sense of~(\ref{domain-convergence}) as $k \to \infty$
($\delta_k \to 0$)
and is essentially the same as the continuum one but with exploration paths
instead of the $\text{SLE}_6$ curves.

\begin{figure}[!ht]
\begin{center}
\includegraphics[width=8cm]{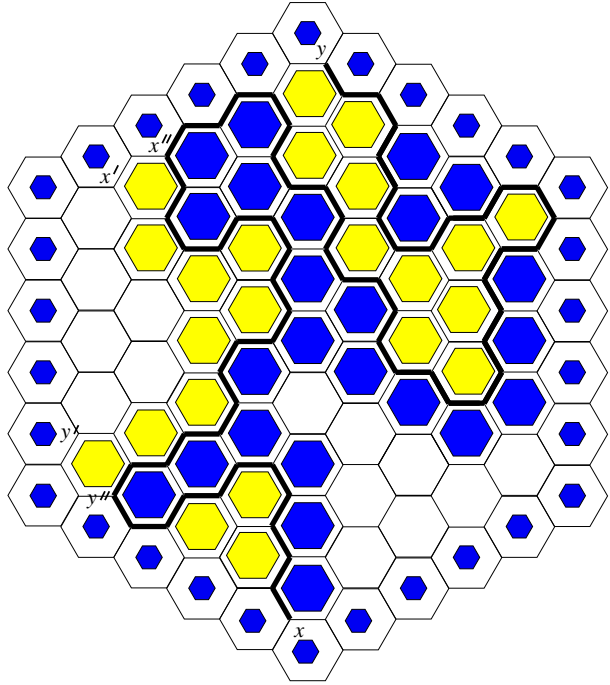}
\caption{First step of the construction of the outer contour
of a cluster of yellow (light in the figure) hexagons
consisting of an exploration (heavy line) from $x$ to $y$.
The outer layer of hexagons does not belong to the domain
where the explorations are carried out, but represents its
monochromatic blue external boundary.
$x''$ and $y''$ are the ending and starting points of an
``excursion" that determines a new domain $D'$, and $x'$ and $y'$
are the vertices where the edges that separate the yellow and
blue portions of the external boundary of $D'$ intersect $\partial D'$.
The second step will consist of an exploration process in $D'$ from
$x'$ to $y'$.}
\label{fig3-sec4}
\end{center}
\end{figure}

This raises the question of how to define an exploration process and obtain
an exploration path in a lattice-domain with monochromatic boundary conditions.
The basic idea is that away from the boundary, the exploration process does not
``know" the boundary conditions.
For two given points $x$ and $y$ on the boundary of a lattice-domain with, say,
blue boundary conditions, split the boundary into two arcs, the counterclockwise
arc $\overline{xy}$ and the the counterclockwise arc $\overline{yx}$.
Then, one can run an exploration process from $x$ to $y$ with the usual rule
inside the domain and on the counterclockwise arc $\overline{xy}$, while
``pretending" that the counterclockwise arc $\overline{yx}$ is colored yellow
(see Fig.~\ref{fig3-sec4}).

If we run such an exploration process in $D_k$ and then look at
the hexagons that have not yet been explored, we will see several
disjoint lattice subdomains, all of which are lattice-Jordan.
This amounts to removing the ``fattened" exploration path consisting of
the exploration path
$\gamma_k \equiv \gamma^{\delta_k}_{D_k,x,y}$
itself and the hexagons immediately to its right and to its left.

The resulting lattice-Jordan subdomains are of four types, which may be
usefully thought of in terms of their external boundaries:
(1) those components  whose boundary contains both sites in the fattened
exploration path and in $\partial_{yx}^k$, the counterclockwise portion
between $y$ and $x$ of the boundary of $D_k$,
(2) the analogous components with $\partial_{yx}^k$ replaced by the other
boundary portion $\partial_{xy}^k$,
(3) those components whose boundary only contains yellow hexagons from
the fattened exploration path and finally
(4) the analogous components whose boundary only contains blue hexagons
from the fattened exploration path.

Notice that the components of type~1 are the only ones with mixed (partly
blue and partly yellow) boundary conditions, while all other components
have monochromatic (blue or yellow) boundary conditions; type~1 components
are special because we have taken blue boundary conditions on $D_k$ while
the exploration path has yellow on its left and blue on its right.
Because of the mixed boundary conditions, each lattice subdomain of type~1 must
contain an interface between the two boundary points where the color changes.
It is also clear that to find such an interface one has to start an exploration
process at one of the two boundary points where the color changes (the two
choices give the same exploration path).

\begin{figure}[!ht]
\begin{center}
\includegraphics[width=8cm]{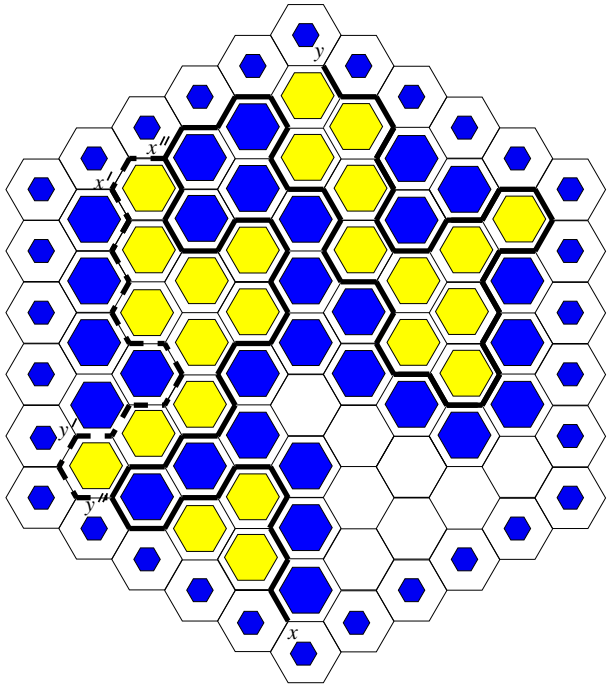}
\caption{Second step of the construction of the outer contour
of a cluster of yellow (light in the figure) hexagons consisting
of an exploration from $x'$ to $y'$ whose resulting path (heavy
broken line) is pasted to  
a portion of the previous exploration path with the help of
the edges (indicated again by a heavy broken line) 
between $x'$ and $x''$ and between $y'$ and $y''$ in such a way
as to obtain a loop around a yellow cluster (light in the figure)
touching the boundary portion $\partial_{yx}^k$.}
\label{fig4-sec4}
\end{center}
\end{figure}

If we run such an exploration process inside a lattice subdomain $D'_k$
of type~1 and paste it to a portion of $\gamma_k$ as in Fig.~\ref{fig4-sec4},
we obtain a loop corresponding to the interface surrounding a yellow cluster
that touches $\partial_{yx}^k$.
If we then again remove the fattened
exploration path, $D'_k$ is split into various components, but this
time those lattice subdomains all have monochromatic boundary conditions.

If we do the same in each subdomain of type~1,
we obtain a collection of loops.
Moreover, all the lattice subdomains of $D_k$
of non-explored hexagons then have monochromatic boundary conditions.
Thus we can iterate the whole procedure inside each of those
lattice subdomains, until we have found all the interfaces contained in $D_k$.

The similarity between this construction and the continuum one
of the CNL process should be apparent.
To continue the proof one needs first to show that the
exploration paths used in the lattice construction converge to chordal
$\text{SLE}_6$ curves.
The first step is a simple application of Theorem~\ref{thm-conv-to-sle}
to the first exploration path
$\gamma_k=\gamma^{\delta_k}_{D_k,x_k,y_k}$, where
$D_k,x_k,y_k$ are chosen so that $D_k$ converges to $D$ and $x_k$ and $y_k$
converge to the $a$ and $b$ of the continuum construction.
However, in order to iterate this step and apply Theorem~\ref{thm-conv-to-sle}
again, we need to also show that the subdomains of the lattice
construction converge to those of the continuum construction.

The convergence in distribution of $\gamma_k$ to $\gamma=\gamma_{D,a,b}$
implies that we can find versions of $\gamma_k$ and $\gamma$ on some
probability space $(\Omega,{\cal B},{\mathbb P})$ such that
$\gamma_k(\omega)$ converges to $\gamma(\omega)$ for all
$\omega \in \Omega$.
Using the coupling, $\gamma_k$ and $\gamma$, for $\delta_k$ small,
are close in the sense of~(\ref{distance}).
This is, however, not sufficient. If we want to conclude convergence
of the subdomains, we need that wherever $\gamma$ touches the boundary
of $D$, $\gamma_k$ touches the boundary of $D_k$ nearby.
Closeness in the sense of~(\ref{distance}) does not ensure this but
only that $\gamma_k$ gets close to the boundary $\partial D_k$.

Note that, if $\gamma_k$ gets within distance $R_1$ of some point $z$
on $\partial D_k$ without touching $\partial D_k$ within distance
$R_2$ of $z$, with $R_2>R_1>\delta_k$, considering the fattened version
of $\gamma_k$ shows the existence of two $(\delta_k{\cal T})$-paths of
one color, say yellow, and one $(\delta_k{\cal T})$-path of the other
color, blue, crossing the annulus of inner radius $R_1$ and outer radius
$R_2$ centered at $z$.

In~\cite{cn1}, where the construction of the CN loops is carried out in
the unit disc $\mathbb D$, the problem is solved by using the fact that
$\mathbb D$ is convex and resorting to an upper bound (see, e.g.,~\cite{lsw5})
on the probability that three disjoint monochromatic $\cal T$-paths cross
a semi-annulus in a half-plane.
The bound shows that such ``three arm" events do not occur in the scaling
limit $\delta \to 0$, implying that, as $k \to \infty$ ($\delta_k \to 0$),
the ($\limsup$ of the) probability that $\gamma_k$ gets within distance
$R_1$ of {\it any} $z \in \partial D_k$ without touching the boundary within
distance $R_2$ of $z$ goes to zero as $R_1 \to 0$
for all (fixed) $R_2 > 0$.

We cannot use that bound here, since $D$ is not necessarily convex
(and even if it were, the $D'$ domains of Theorems~\ref{thm-conformal-invariance}
and ~\ref{thm-cle} will not generally be convex).
Instead, we will use the continuity of Cardy's formula with respect to small
changes in the shape of the domain.
We postpone this issue until later and proceed with the sketch of the proof
assuming that $\gamma_k$ does not get close to the boundary of the domain
without touching it nearby (probably).

Then the boundaries of the lattice/continuum subdomains
obtained after running the first (coupled) exploration path/$\text{SLE}_6$
curve are close to each other in the metric~(\ref{distance}).
I.e., we can match lattice and continuum subdomains, at least
for those whose diameter is larger than some $\varepsilon_k$
which depends on $\delta_k$.
It is important that, as $k \to \infty$ (and $\delta_k \to 0$),
we can let $\varepsilon_k \to 0$.

If we run an exploration process inside a (large) lattice subdomain
$D'_k$ converging to a continuum subdomain $D'$, Theorem~\ref{thm-conv-to-sle}
allows us to conclude that the exploration path $\gamma'_k$ in
$D'_k$ converges to the $\text{SLE}_6$ curve $\gamma'$ in $D'$
from $a'$ to $b'$, provided that the starting and ending points $x'_k$
and $y'_k$ of the exploration process are chosen so that they converge
to $a'$ and $b'$ respectively as $k \to \infty$.
We can now work with coupled versions of $\gamma'_k$ and $\gamma'$
and repeat the above argument with the new subdomains that they
produce, obtaining again a match (with high probability).

This allows us to keep the lattice and continuum constructions coupled,
which ensures in particular that the $(\delta_k{\cal T})$-loops obtained
in the lattice construction converge, as $\delta_k \to 0$, to the
loops obtained in the continuum construction.

For any fixed $\delta_k$, it is clear that the lattice construction
eventually finds all the boundary loops.
However, to conclude that the CNL process is indeed the scaling limit
of the collection of all interfaces, we need to show that, for any
$\varepsilon>0$, the number of steps of the discrete construction
needed to find all the loops of diameter at least $\varepsilon$ does
not diverge as $k \to \infty$ (otherwise some loops would never be
found in the scaling limit).

In~\cite{cn1}, this is resolved using percolation arguments (that make
use of the RSW theorem~\cite{russo, sewe} and FKG inequalities) to show
that the size of the subdomains has a bounded away from zero probability
of decreasing significantly at each iteration.
We point out that the argument used in~\cite{cn1}, where the construction
of the CN loops in carried out in the unit disc, is independent of the
actual shape of the domain so that
it can be applied to the present situation.
Since that argument is long, we will not repeat it here.

\begin{figure}[!ht]
\begin{center}
\includegraphics[width=7cm]{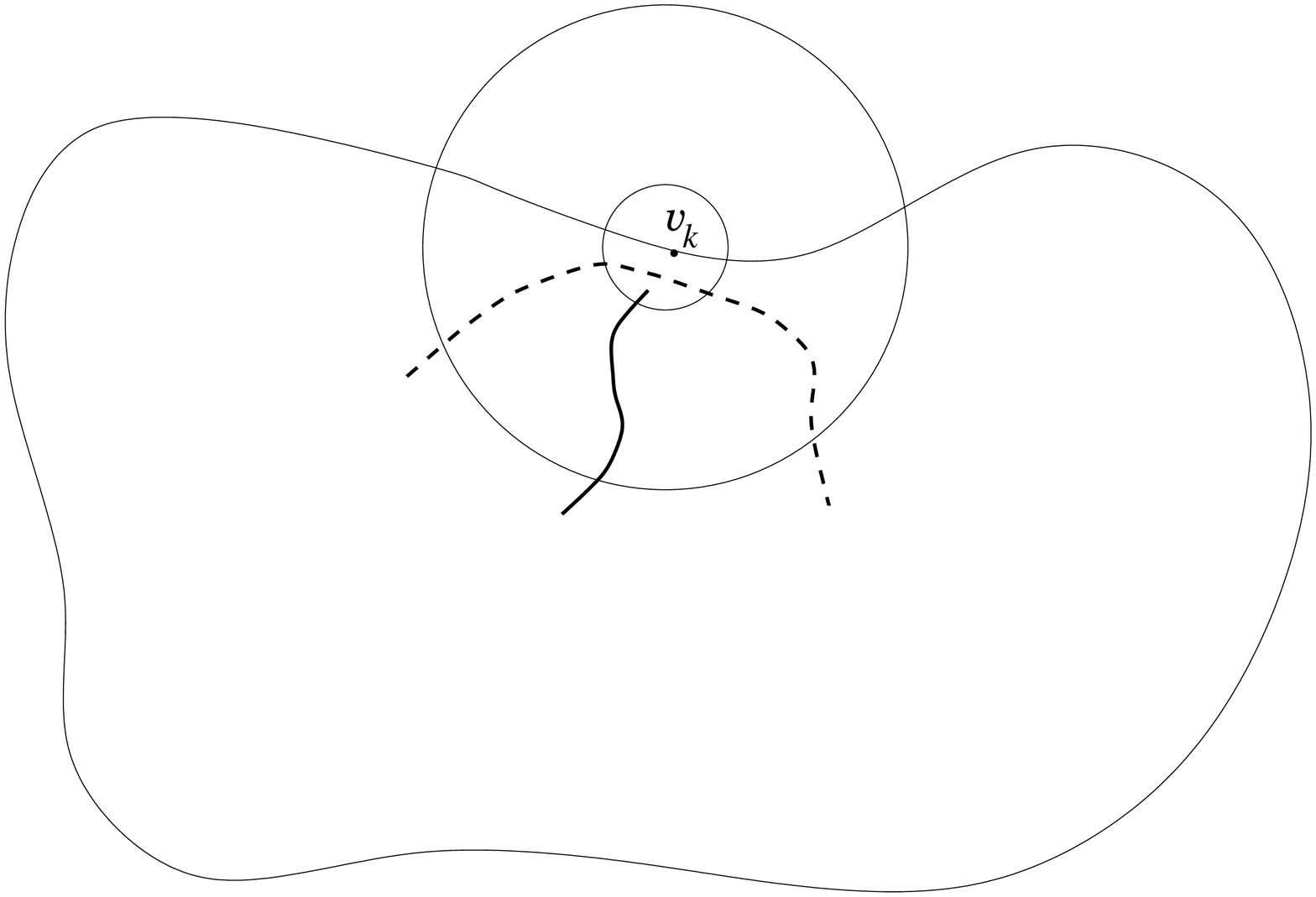}
\caption{The figure shows a blue $(\delta_k{\cal T})$-path
(heavy full line) crossing the partial annulus
$D_k \cap \{ B(v_k,R) \setminus B(v_k,r) \}$ that fails to
connect to $\partial D_k$ near $v_k$ because it is blocked by
a yellow $(\delta_k{\cal T})$-path (heavy dashed line) that
twice crosses the annulus $B(v_k,R) \setminus B(v_k,r)$.}
\label{fig-mushroom}
\end{center}
\end{figure}

Returning to the problem of ``close encounters" of $\gamma_k$
with $\partial D_k$, we will try to provide the intuition on which
the proof of touching is
based.
Suppose, by contradiction, that $\gamma_k$ enters the disc
$B(v_k,\varepsilon_k)$ of radius $\varepsilon_k$ centered at
$v_k \in \partial D_k$ without touching $\partial D_k$ inside
the disc $B(v_k,r)$ of radius $r$, and that $\varepsilon_k \to 0$.
As $k \to \infty$, $D_k \to D$ and we can assume by compactness
that $v_k$ converges to some $v \in \partial D$.
Considering the fattened version of $\gamma_k$ shows the existence
of two $(\delta_k{\cal T})$-paths of one color, say yellow, and one
$(\delta_k{\cal T})$-path of the other color, blue, crossing the
annulus $B(v_k,r) \setminus B(v_k,\varepsilon_k)$ (see Fig.~\ref{fig-mushroom}).

Assume for simplicity that $v$ is far enough from $a$ and
$b$ so that $a,b \notin \overline{B(v,R)}$ for some $R>r$, and
consequently $x_k,y_k \notin B(v_k,R)$ for $k$ large enough.
Then, in the domain
$D_k \cap \{ B(v_k,R) \setminus B(v_k,r) \}$ there is a blue crossing
between a certain portion $J_k$ of the circle of radius $R$ centered
at $v_k$ and a certain portion $J'_k$ of the circle of radius $r$
centered at $v_k$.
If we consider instead the domain $D_k \cap B(v_k,R)$, there is no blue
crossing between $J_k$ and the portion of $\partial D_k \cap B(v_k,r)$
containing $v_k$ (see Fig.~\ref{fig-mushroom}).
If this discrepancy persists as $k \to \infty$, it must show up
in the scaling limit of crossing probabilities for the domains
$D \cap \{ B(v,R) \setminus B(v,r) \}$ and $D \cap B(v,R)$.
On the other hand, since $\varepsilon_k \to 0$, we can take $r$ very
small, and so $D \cap \{ B(v,R) \setminus B(v,r) \}$ is very close
to $D \cap B(v,R)$ so that the ``crossing probabilities" in the two
domains between the corresponding arcs, given in the continuum by
Cardy's formula, should be very close.
This follows from the continuity of Cardy's formula with respect to
the shape of the domain and the positions of the boundary arcs (see,
e.g., Lemma A.2 of~\cite{cn2}).

Using this idea, one can show that the assumption that $\gamma_k$ comes
close to $\partial D_k$ without touching it nearby produces a contradiction.
Although the idea outlined above is relatively simple, the arguments needed
to obtain a contradiction are rather involved (see Lemmas 7.1, 7.2, 7.3
and 7.4 of~\cite{cn2}), so we will not present them here, except for a brief
discussion following Lemma~\ref{combined-boundaries} below. \\

\noindent {\bf Sketch of the Proof of Theorem~\ref{thm-conformal-invariance}.}
In order to prove the claim, we will define a lattice
construction inside $D'$ coupled to the continuum
construction inside $D$, by means of the conformal map
$f$ from $D$ to $D'$.
Roughly speaking, this new lattice construction for $D'$
is one in which the $(x,y)$ pairs at each step are chosen
to be close to the $(f(a),f(b))$ points in $D'$ mapped
from $D$ via $f$, where the pairs $(a,b)$ are those that
appear at the corresponding steps of the continuum
construction inside $D$.

More precisely, let $\gamma_{(1)}$ be the first $\text{SLE}_6$
curve in $D$ from $a_{(1)}$ to $b_{(1)}$.
Because of the conformal invariance of $\text{SLE}_6$,
the image $f(\gamma_{(1)})$ of $\gamma_{(1)}$ under $f$ is a
curve
distributed as the trace of chordal $\text{SLE}_6$ in $D'$
from $f(a_{(1)})$ to $f(b_{(1)})$.
Therefore, the exploration path $\gamma^{\delta}_{(1)}$ inside
$D'$ from $x_{(1)}$ to $y_{(1)}$, chosen so that they converge to
$f(a_{(1)})$ and $f(b_{(1)})$ respectively as $\delta \to 0$,
converges in distribution to $f(\gamma_{(1)})$, as $\delta \to 0$,
which means that there exists a coupling between $\gamma^{\delta}_{(1)}$
and $f(\gamma_{(1)})$ such that the curves stay close for $\delta$ small.

We see that one can use the same strategy as in the
sketch of the proof of Theorem~\ref{thm-conv-to-sle}, and
obtain a lattice construction whose exploration paths are
coupled to the $\text{SLE}_6$ curves in $D'$ that are the
images under $f$ of the $\text{SLE}_6$ curves in $D$.
Then, for this discrete construction, the scaling limits
of the exploration paths will be distributed as the images
of the $\text{SLE}_6$ curves in $D$.

To conclude the proof, we should show that
the lattice construction inside $D'$ defined above finds
all the boundaries in a number of steps that is bounded in
probability as $\delta \to 0$. But this is essentially equivalent to
the analogous claim in the sketch
of the proof of Theorem~\ref{thm-conv-to-sle}.
Thus the scaling limit, as $\delta \to 0$, of this
new lattice construction for $D'$ gives the CNL process in
$D'$, which by construction is distributed like the image
under $f$ of the CNL process in $D$. \\

\noindent {\bf Sketch of the Proof of Theorem~\ref{thm-cle}.}
Let $a,b \in \partial D$ be the endpoints of $\Gamma$ in clockwise order,
i.e., $\Gamma=\overline{ba}$ is the counterclockwise arc of $\partial D$
from $b$ to $a$.
As explained at the end of Sec.~\ref{sec-full}, the random curve $\gamma$ obtained
by pasting together the maximal excursions $\overline{xy}(L)$ from $\overline{ba}$,
for $L \in {\cal L}_{D,\Gamma}$, is distributed like chordal $\text{SLE}_6$ in $D$
from $a$ to $b$.
Indeed, removing $\gamma$ from $D$ is equivalent (in distribution) to the first
step of the algorithmic construction presented in Sec.~\ref{sec-full} to produce a
realization of the CNL process, if we choose $a$ and $b$ with $\overline{ba}=\Gamma$
as starting and ending points of the first $\text{SLE}_6$ curve of the construction.

Note that $\gamma$ is in
${\cal L}_{D,\Gamma}^* \equiv \overline{\cup \{ L:L \in {\cal L}_{D,\Gamma} \}}$,
and the remaining pieces of ${\cal L}_{D,\Gamma}^*$
are all in (the closures of) subdomains of $D \setminus \gamma$ of type~1.
If we condition on $\gamma$ and run the algorithmic construction described in
Sec.~\ref{sec-full} inside a subdomain of $D \setminus \gamma$ of type~2, 3 or 4,
we get an independent CNL process or, by Theorem~\ref{thm-conformal-invariance}, an
independent copy of ${\cal L}_D$ conformally mapped to that domain.
This already proves part of the claim.

Consider now a subdomain $D'$ of $D \setminus \gamma$ of type~1 and let
$a',b'$ be the endpoints of the excursion that generated $D'$.
Part of $\partial D'$ is in $\partial D$ and we choose $a', b'$ so that
the counterclockwise arc $\Gamma'=\overline{b'a'} \subset \partial D$
is that part of $\partial D'$.
The excursion that generated $D'$ is part of a loop $L'$ whose other
``half" is in $D'$ and runs from $b'$ to $a'$.
We know from the construction of Sec.~\ref{sec-full} that if we trace the
``half" of $L'$ contained in $D'$ from $b'$ to $a'$ we get a curve $\gamma'$
distributed like chordal $\text{SLE}_6$ in $D'$ from $b'$ to $a'$.
Note that $\gamma'$ is contained in ${\cal L}_{D,\Gamma}^*$.

The subdomains of $D' \setminus \gamma'$ are of two types: (I) those whose
boundary does not contain a portion of $\partial D$ and (II) those whose
boundary does contain a portion, $\Gamma''=\overline{b''a''} \subset \partial D$,
of $\Gamma$.
If we condition on $\gamma$ and $\gamma'$ and run the algorithmic construction
described in Sec.~\ref{sec-full} inside a subdomain of $D' \setminus \gamma'$
of type~I, we get an independent CNL process or, by Theorem~\ref{thm-conformal-invariance},
an independent copy of ${\cal L}_D$ conformally mapped to that domain.

The remaining pieces of ${\cal L}_{D,\Gamma}^*$
are all contained inside the (closures of) domains of type~II (for all the
subdomains of $D \setminus \gamma$ of type~1).
Inside each subdomain $D''$ of type~II, the CN loops that touch $\Gamma''$
are contained in ${\cal L}_{D,\Gamma}^*$
and can be used to obtain a curve $\gamma''$ distributed like chordal
$\text{SLE}_6$ in $D''$ from $a''$ to $b''$ by pasting together maximal
excursions as above (and at the end of Sec.~\ref{sec-full}).
It should now be clear how to complete the argument by iterating the steps
described above inside each subdomain $D''$.


\section{Convergence of Exploration Path to $\text{SLE}_{\bf 6}$} \label{sec-convergence-sle}



We begin discussing the proof of Theorem~\ref{thm-conv-to-sle} by noting
(like in the proof of Theorem~\ref{thm-full-scaling-limit} discussed in
Sec.~\ref{sec-convergence-full}) that from the work of Aizenman-Burchard~\cite{ab},
the family of distributions of $\gamma^{\delta}_{D,a,b}$  is tight (as
$\delta \to 0$, in the sense of the metric~(\ref{distance})) and so there
is convergence along subsequences $\delta_k \to 0$.
We write, in simplified notation, $\gamma_k \to \tilde{\gamma}$ along such
a convergent subsequence.
What needs to be proved is that the distribution $\tilde{\mu}$ of $\tilde{\gamma}$
is that of $\gamma^{\text{SLE}_6}$, the trace of chordal $\text{SLE}_6$ in $D$
from $a$ to $b$.

We next discuss how much information about $\tilde{\mu}$ can be extracted
from Cardy's formula for crossing probabilities.
We note that there are versions of Smirnov's result on convergence of crossing
probabilities to Cardy's formula that allow the domains being crossed and the
target boundary arcs to vary as $\delta \to 0$. Theorem~3 of~\cite{cn2} is such
a version that suffices for our purposes.
Let $D_t \equiv D \setminus \tilde K_t$ denote the (unique) connected
component of $D \setminus \tilde\gamma[0,t]$ whose closure contains $b$,
where $\tilde K_t$, the {\bf filling} of $\tilde\gamma[0,t]$, is a closed
connected subset of $\overline D$.
$\tilde K_t$ is called a {\bf hull} if it satisfies the condition
\begin{equation} \label{hulls}
\overline{\tilde K_t \cap D} = \tilde K_t.
\end{equation}
We will consider curves $\tilde\gamma$ such that $\tilde K_t$ is a hull for
each $t$, although here we only consider $\tilde K_T$ at certain
stopping times $T$.

\begin{figure}[!ht]
\begin{center}
\includegraphics[width=7cm]{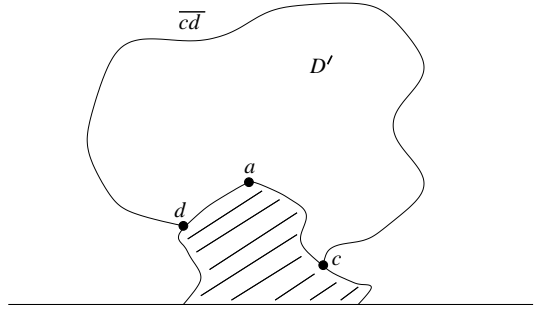}
\caption{$D$ is the upper half-plane $\mathbb H$ with the shaded portion removed,
$b=\infty$, $C'$ is an unbounded subdomain, and $D' = D \setminus C'$ is indicated
in the figure. The counterclockwise arc $\overline{cd}$ indicated in the figure
belongs to $\partial D'$.}
\label{fig1-sec7}
\end{center}
\end{figure}

Let $C' \subset D$ be a closed subset of $\overline D$ such that $a \notin C'$,
$b \in C'$, and $D' = D \setminus C'$ is a bounded simply connected domain
whose boundary contains the counterclockwise arc $\overline {cd}$ that does
not belong to $\partial D$ (except for its endpoints $c$ and $d$ -- see
Fig.~\ref{fig1-sec7}).
Let $T'=\inf \{ t:\tilde K_t \cap C' \neq \emptyset \}$ be the first time
that $\tilde\gamma(t)$ hits $C'$ and assume that the filling $\tilde K_{T'}$
of $\tilde\gamma[0,T']$ is a hull.
We say that the hitting distribution of $\tilde\gamma(t)$
at the stopping time $T'$ is determined by Cardy's
formula (see~(\ref{cardy-formula})) if, for any $C'$ and any counterclockwise arc
$\overline{xy}$ of $\overline{cd}$, the probability that $\tilde\gamma$ hits $C'$
at time $T'$ on $\overline{xy}$ is given by
\begin{equation}
{\mathbb P}(\tilde\gamma(T') \in \overline{xy}) = \Phi_{D'}(a,c;x,d) - \Phi_{D'}(a,c;y,d).
\end{equation}

\begin{figure}[!ht]
\begin{center}
\includegraphics[width=7cm]{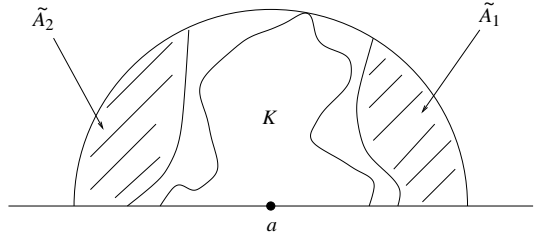}
\caption{Example of a hull $K$ and a set $\tilde A_1 \cup \tilde A_2$
(shaded regions) in $\cal A$.
Here, $D = {\mathbb H}$ and $D'$ is the semi-disc centered at $a$.}
\label{fig2-sec7}
\end{center}
\end{figure}

We want to relate the distribution of $\tilde K_{T'}$ to the distribution
of hitting \emph{locations} for a family of $C''$'s related to $C'$.
To explain, consider the set $\tilde{\cal A}$ of closed subsets $\tilde A$
of $\overline{D'}$ that do not contain $a$ and such that
$\partial \tilde A \setminus \partial D'$ is a simple (continuous) curve
contained in $D'$
except for its endpoints, one of which is on $\partial D' \cap D$ and the other is
on $\partial D$ (see Fig.~\ref{fig2-sec7}).
Let $\cal A$ be the set of closed subsets of $\overline{D'}$ of the form
$\tilde A_1 \cup \tilde A_2$, where $\tilde A_1, \tilde A_2 \in \tilde{\cal A}$
and $\tilde A_1 \cap \tilde A_2 = \emptyset$.

It is easy to see that if the hitting distribution of
$\tilde\gamma(T')$ is
determined by Cardy's formula, then the probabilities of events of the form
$\{\tilde K_{T'} \cap A = \emptyset\}$ for $A \in {\cal A}$
are also determined by Cardy's formula in the following way.
Let $A \in {\cal A}$ be the union of $\tilde A_1, \tilde A_2 \in \tilde{\cal A}$,
with $\partial \tilde A_1 \setminus \partial D'$ given by a curve from
$u_1 \in \partial D' \cap D$ to $v_1 \in \partial D$ and
$\partial \tilde A_2 \setminus \partial D'$ given by a curve from
$u_2 \in \partial D' \cap D$ to $v_2 \in \partial D$; then, assuming that
$a$, $v_1$, $u_1$, $u_2$, $v_2$ are ordered counterclockwise around $\partial D'$,
\begin{equation}
{\mathbb P}(\tilde K_{T'} \cap A = \emptyset) =
\Phi_{D' \setminus A}(a,v_1;u_1,v_2,) - \Phi_{D' \setminus A}(a,v_1;u_2,v_2).
\end{equation}
The probabilities of such events determine uniquely the distribution of the hull
(for more detail, see Sec.~5 of~\cite{cn2}).
Thus we have the following useful lemma, since the hitting distribution for
$\text{SLE}_6$ is determined by Cardy's formula~\cite{lsw1}.

\begin{lemma} \label{hull}
If $\tilde K_{T'}$ is a hull and the hitting distribution of $\tilde\gamma$
at the stopping time $T'$
is determined by Cardy's formula, then $\tilde K_{T'}$ is distributed like
the corresponding hull of $\gamma^{\text{\emph{SLE}}_6}$.
\end{lemma}

We next define the sequence of hitting times for $\tilde{\gamma}$ that will be used
to compare it to $\gamma^{\text{SLE}_6}$. They involve conformal maps of semi-balls
(i.e., half-disks) in the upper half-plane.
Let $\tilde f_0$ be a conformal map from the upper half-plane $\mathbb H$
to $D$ such that ${\tilde f_0}^{-1}(a)=0$ and ${\tilde f_0}^{-1}(b)=\infty$.
(Since $\partial D$ is a continuous curve, the map ${\tilde f_0}^{-1}$ has a
continuous extension from $D$ to $D \cup \partial D$ 
and, by a slight abuse of notation, we do not
distinguish between ${\tilde f_0}^{-1}$ and its extension; the same applies
to $\tilde f_0$.)
These two conditions determine $\tilde f_0$ only up to a scaling factor.
For $\varepsilon>0$ fixed, let
$C(u,\varepsilon) = \{ z:|u-z|<\varepsilon \} \cap {\mathbb H}$
denote the semi-ball of radius $\varepsilon$ centered at $u$ on
the real line and let $\tilde T_1=\tilde T_1(\varepsilon)$ denote
the first time $\tilde\gamma(t)$ hits $D \setminus \tilde G_1$,
where $\tilde G_1 \equiv \tilde f_0(C(0,\varepsilon))$.
Define recursively $\tilde T_{j+1}$ as the first time
$\tilde\gamma[\tilde T_j,\infty)$ hits $\tilde D_{\tilde T_j} \setminus \tilde G_{j+1}$,
where $\tilde D_{\tilde T_j} \equiv D \setminus \tilde K_{\tilde T_j}$,
$\tilde G_{j+1} \equiv \tilde f_{\tilde T_j}(C(0,\varepsilon))$, and
$\tilde f_{\tilde T_j}$ is a conformal map from $\mathbb H$ to $\tilde D_{\tilde T_j}$
whose inverse maps $\tilde\gamma(\tilde T_j)$ to $0$ and $b$ to $\infty$.
We also define $\tilde\tau_{j+1} \equiv \tilde T_{j+1} - \tilde T_j$,
so that $\tilde T_j = \tilde\tau_1+\ldots+\tilde\tau_j$.
We choose $\tilde f_{\tilde T_j}$ so that its inverse is the composition
of the restriction of $\tilde{f_0}^{-1}$ to $\tilde D_{\tilde T_j}$ with
$\tilde \varphi_{\tilde T_j}$, where $\tilde \varphi_{\tilde T_j}$ is the
unique conformal transformation from
${\mathbb H} \setminus \tilde{f_0}^{-1}(\tilde K_{\tilde T_j})$ to $\mathbb H$
that maps $\infty$ to $\infty$ and $\tilde{f_0}^{-1}(\tilde\gamma(\tilde T_j))$
to the origin of the real axis, and has derivative at $\infty$ equal to $1$.

Notice that $\tilde G_{j+1}$ is a bounded simply connected domain chosen
so that the conformal transformation which maps $\tilde D_{\tilde T_j}$
to $\mathbb H$ maps $\tilde G_{j+1}$ to the semi-ball $C(0,\varepsilon)$
centered at the origin on the real line. With these definitions, we consider
the (discrete-time) stochastic process
$\tilde X_j \equiv (\tilde K_{\tilde T_j}, \tilde\gamma(\tilde T_j))$
for $j=1,2,\ldots$ .
Analogous quantities can be defined for the trace of chordal $\text{SLE}_6$.
They are indicated by the superscript $\text{SLE}_6$;
we choose $f_0^{\text{SLE}_6}=\tilde f_0$, so that $G_1^{\text{SLE}_6}=\tilde G_1$.
Our aim is to prove that the variables $\tilde X_1, \tilde X_2,\dots$
are (jointly) equidistributed with the corresponding $\text{SLE}_6$ hull and
tip variables $X_1^{\text{SLE}_6}, X_2^{\text{SLE}_6},\dots$ .
By letting $\varepsilon \to 0$, this will directly yield that $\tilde{\gamma}$
is equidistributed with $\gamma^{\text{SLE}_6}$ as desired. Since $\gamma_k$ converges
in distribution to $\tilde\gamma$, we can find coupled versions of $\gamma_k$
and $\tilde\gamma$ on some probability space $(\Omega,{\cal B},{\mathbb P})$
such that $\gamma_k$ converges to $\tilde\gamma$ for all $\omega \in \Omega$;
in the rest of the proof we work with these new versions which, with a slight
abuse of notation, we denote with the same names as the original ones.

For each $k$, let $K^k_t$ denote the filling (or {\bf lattice hull}) at time
$t$ of $\gamma_k$, i.e., the set of hexagons that at time $t$ have been explored
or have been disconnected from $b$ by the exploration path.
Let now $f^k_0$ be a conformal transformation that maps $\mathbb H$ to
$D_k \equiv D^{\delta_k}$ such that $(f^k_0)^{-1}(a_k)=0$ and
$(f^k_0)^{-1}(b_k) = \infty$ and let $T^k_1=T^k_1(\varepsilon)$ denote the first
exit time of $\gamma^{\delta_k}_k(t)$ from $G^k_1 \equiv f_0^k(C(0,\varepsilon))$
defined as the first time
that $\gamma_k$ intersects the image under $f_0^k$
of the semi-circle $\{ z : |z| = \varepsilon \} \cap {\mathbb H}$.
Define recursively $T^k_{j+1}$ as the first exit time of $\gamma^{\delta_k}_k[T^k_j,\infty)$
from $G^k_{j+1} \equiv f^k_{T^k_j}(C(0,\varepsilon))$, where $f^k_{T^k_j}$
is a conformal map from $\mathbb H$ to $D_k \setminus K^k_{T^k_j}$ whose
inverse maps $\gamma_k(T^k_j)$ to $0$ and $b_k$ to $\infty$.
The maps
$f^k_{T^k_j}$, for $j \geq 1$, are defined only up to a scaling factor.
We also define $\tau^k_{j+1} \equiv T^k_{j+1} - T^k_j$, so that
$T^k_j=\tau^k_1+\ldots+\tau^k_j$, and the (discrete-time) stochastic process
$X^k_j \equiv (K^k_{T^k_j},\gamma^{\delta_k}_k(T^k_j))$ for $j=1,2,\ldots$ .

We want to show recursively that, for any $j$, as $k \to \infty$,
$\{ X^k_1,\ldots,X^k_j \}$ converge jointly in distribution to
$\{ \tilde X_1,\ldots,\tilde X_j \}$.
By recursively applying convergence of crossing probabilities
to Cardy's formula (i.e., Theorem~3 of~\cite{cn2}) and Lemma~\ref{hull}, we
will then be able to conclude, as explained in more detail below, that
$\{ \tilde X_1, \tilde X_2, \ldots \}$ are jointly equidistributed with
the corresponding $\text{SLE}_6$ hull variables (at the corresponding stopping times)
$\{ X_1^{\text{SLE}_6},X_2^{\text{SLE}_6},\ldots \}$.

The zeroth step consists in noticing that the convergence of $(D_k,a_k,b_k)$
to $(D,a,b)$ as $k \to \infty$ allows us to
select a sequence of conformal maps $f^k_0$ that
converge to $f_0^{\text{SLE}_6}=\tilde f_0$ uniformly in $\overline{\mathbb H}$
as $k \to \infty$, which implies that the boundary $\partial G^k_1$
of $G^k_1=f^k_0(C(0,\varepsilon))$ converges to the boundary $\partial\tilde G_1$
of $\tilde G_1 = \tilde f_0(C(0,\varepsilon))$ in the uniform metric on continuous
curves (see Corollary~A.2 of~\cite{cn2}).

The next lemma is the technical heart of the proof. It basically allows us
to interchange the scaling limit $\delta \to 0$ and the process of filling
(which generates hulls) by declaring that the hull of the limiting curve is
the limit of the (lattice) hulls.
The proof of the lemma involves extensive use of nontrivial results from
percolation theory.
Although the lemma is stated here in the framework of the first step of
the proof where we are analyzing convergence of $X_1^k$ to $\tilde X_1$,
essentially the same lemma can be applied sequentially to the convergence
of $X_j^k$ conditioned on $\{ X^k_1,\ldots,X^k_{j-1} \}$.

\begin{lemma} \label{combined-boundaries}
$(\gamma_k, K^k_{T^k_1})$ converges in distribution to
$(\tilde \gamma, \tilde K_{\tilde T_1})$ as $k \to \infty$.
Furthermore $\tilde K_{\tilde T_1}$ is a.s. a hull
equidistributed with the hull $K^{\text{\emph{SLE}}_6}_{T_1}$ of
$\text{\emph{SLE}}_6$ at the corresponding stopping time $T_1$.
\end{lemma}

Proving the first claim,
that for the exploration path $\gamma_k$ in $G^k_1$ one can
interchange the limit $k \to \infty$ (${\delta}_k \to 0$) with the
process of filling, requires showing two things about the exploration
path: (1) the return of a (macroscopic) segment of the path close to
an earlier segment (and away from $\partial G^k_1$) without nearby
(microscopic) touching does not occur (probably), and (2) the close approach of a
(macroscopic) segment of the path to $\partial G^k_1$ without nearby (microscopic)
touching either of $\partial G^k_1$ itself or else of another segment of the path
that touches $\partial G^k_1$ does not occur (probably).
If $G^k_1$ (or more accurately, its limit $\tilde G_1$) were
replaced by a convex domain like the unit disk, these could be
controlled by known estimates on probabilities of six-arm events in
the full plane for (1) and of three-arm events in the half-plane for (2).
But $\tilde G_1$ is not in general convex and then
the three-arm event argument for (2) appears to break down. The replacement
in~\cite{cn2} is the use of several lemmas in Sec.~7 there.
Basically, these control (2) by a novel argument about ``mushroom events''
in $\tilde G_1$, 
which is based on continuity of Cardy's formula with respect to changes
in $\partial \tilde G_1$.
Roughly speaking, mushroom events are ones where (in the limit $k \to \infty$)
there is a macroscopic monochromatic path in $\tilde G_1$ just reaching to
$\partial \tilde G_1$, but blocked from it by a macroscopic path in
$\tilde G_1$ of the other color
(see Fig.~\ref{fig-mushroom}).
It is shown in~\cite{cn2} (see Lemma~7.4 there) that mushroom events cannot
occur with positive probability while on the other hand they would occur if
(2) were not the case.
%
The second claim of Lemma~\ref{combined-boundaries} now follows
from Smirnov's result on convergence to
Cardy's formula~\cite{smirnov, smirnov-long} (see also Theorem~3 of~\cite{cn2})
and Lemma~\ref{hull}.

Using Lemma~\ref{combined-boundaries}, the first step of our recursion argument is organized
as follows, where all limits and equalities are in distribution:
\begin{itemize}
\item[(i)] $K^k_{T^k_1} \to \tilde K_{\tilde T_1} = K_{T_1}^{\text{SLE}_6}$ by
Lemma~\ref{combined-boundaries}. 
\item[(ii)] by i), $D_k \setminus K^k_{T^k_1} \to D \setminus \tilde K_{\tilde T_1}
= D \setminus K_{T_1}^{\text{SLE}_6}$.
\item[(iii)] by (ii), $f_{T_1}^{\text{SLE}_6}=\tilde f_{\tilde T_1}$,
and 
we can select
a sequence $f^k_{T^k_1} \to \tilde f_{\tilde T_1} = f_{T_1}^{\text{SLE}_6}$.
\item[(iv)] by (iii), $G^k_2 \to \tilde G_2 = G_2^{\text{SLE}_6}$.
\end{itemize}
At this point, we are in the same situation as at the zeroth step,
but with $G^k_1$, $\tilde G_1$ and $G_1^{\text{SLE}_6}$ replaced by $G^k_2$,
$\tilde G_2$ and $G_2^{\text{SLE}_6}$, and we proceed by induction, as follows.

The next step consists in proving that
$((K^k_{T^k_1},\gamma^{\delta_k}_k(T^k_1)),(K^k_{T^k_2},\gamma^{\delta_k}_k(T^k_2)))$
converges in distribution to
$((\tilde K_{\tilde T_1},\tilde\gamma(\tilde T_1)),(\tilde K_{\tilde T_2},\tilde\gamma(\tilde T_2)))$.
Since we have already proved the convergence of $(K^k_{T^k_1},\gamma^{\delta_k}_k(T^k_1))$
to $(\tilde K_{\tilde T_1},\tilde\gamma(\tilde T_1))$, we claim that all we really need
to prove is the convergence of
$(K^k_{T^k_2} \setminus K^k_{T^k_1},\gamma^{\delta_k}_k(T^k_2))$
to $(\tilde K_{\tilde T_2} \setminus \tilde K_{\tilde T_1},\tilde\gamma(\tilde T_2))$.
To do this, notice that $K^k_{T^k_2} \setminus K^k_{T^k_1}$ is distributed like
the lattice hull of a percolation exploration path inside $D_k \setminus K^k_{T^k_1}$.
Besides, the convergence in distribution of $(K^k_{T^k_1},\gamma^{\delta_k}_k(T^k_1))$
to $(\tilde K_{\tilde T_1},\tilde\gamma(\tilde T_1))$ implies that we can find versions
of $(\gamma^{\delta_k}_k,K^k_{T^k_1})$ and $(\tilde\gamma,\tilde K_{\tilde T_1})$
on some probability space $(\Omega,{\cal B},{\mathbb P})$ such that
$\gamma^{\delta_k}_k(\omega)$ converges to $\tilde\gamma(\omega)$ and
$(K^k_{T^k_1},\gamma^{\delta_k}_k(T^k_1))$ converges to
$(\tilde K_{\tilde T_1},\tilde\gamma(\tilde T_1))$ for all $\omega \in \Omega$.
These two observations imply that, if we work with the coupled versions of
$(\gamma^{\delta_k}_k,K^k_{T^k_1})$ and $(\tilde\gamma,\tilde K_{\tilde T_1})$,
we are in the same situation as before, but with $D_k$ (resp., $D$) replaced by
$D_k \setminus K^k_{T^k_1}$ (resp., $D \setminus \tilde K_{\tilde T_1}$) and $a_k$
(resp., $a$) by $\gamma^{\delta_k}_k(T^k_1)$ (resp., $\tilde\gamma(\tilde T_1)$).
Then, the conclusion that
$(K^k_{T^k_2} \setminus K^k_{T^k_1},\gamma^{\delta_k}_k(T^k_2))$
converges in distribution to
$(\tilde K_{\tilde T_2} \setminus \tilde K_{\tilde T_1},\tilde\gamma(\tilde T_2))$
follows, as before, by arguments like those used for Lemma~\ref{combined-boundaries}.
We can now iterate the above arguments $j$ times, for any $j>1$.
If we keep track at each step of the previous ones, this provides the \emph{joint}
convergence of all the curves and lattice hulls involved at each step.

The proof of Theorem~\ref{thm-conv-to-sle} is concluded by letting $\varepsilon \to 0$.
We note that in this paper we circumvent the use of a ``spatial Markov property"
that played a role in~\cite{cn2} in the $\varepsilon \to 0$ limit.
The point is that that property was proved as a consequence of the equidistribution
of $\tilde X_1, \tilde X_2, \ldots$ with $X^{\text{SLE}_6}_1, X^{\text{SLE}_6}_2, \ldots$
and here we apply the equidistribution directly.
It should be noted however that there needs to be some ``a priori" information
about $\tilde\gamma$ to insure that this equidistribution for each
$\varepsilon>0$ implies equidistribution of $\tilde\gamma$ with $\gamma^{\text{SLE}_6}$.
For example, one could create by hand a process $\hat\gamma$ which behaved like
$\gamma^{\text{SLE}_6}$ except that at random times it retraced back and forth part
of its previous path.
Such a $\hat\gamma$ would have its $\hat X_j$ variables equidistributed with those
of $\text{SLE}_6$ but as a random curve (modulo monotonic reparametrizations) would
not be equidistributed with $\gamma^{\text{SLE}_6}$; it would also not be describable
by a Loewner chain.
Such possibilities can be ruled out by the same arguments as those used in proving
Lemma~\ref{combined-boundaries} (see Lemma~6.4 of~\cite{cn2}). \fbox{} \\

{\bf Acknowledgements.} The authors thank the Kavli Institute for Theoretical
Physics for its hospitality in 2006 when this paper was mostly written.
C.~M.~N. thanks the Clay Mathematics Institute for partial support of his visit
at KITP, the Mathematical Sciences Research Institute, Berkeley for organizing
the 2005 workshop in honor of Henry McKean and he thanks Henry for many
interesting conversations over many years.

\end{document}